\documentclass[11pt,a4paper]{amsart}

\usepackage{amsmath,amssymb,amsthm,mathtools,mathrsfs}
\usepackage[margin=1.1in]{geometry}
\usepackage{microtype}
\usepackage{tikz}
\usetikzlibrary{arrows.meta,positioning}
\usepackage[colorlinks=true,linkcolor=blue,citecolor=blue,urlcolor=blue]{hyperref}

\newtheorem{theorem}{Theorem}[section]

\theoremstyle{definition}

\theoremstyle{remark}

\numberwithin{equation}{section}

\newcommand{\R}{\mathbb R}

\title[LLN for ASM]
{Limiting shape of alternating sign matrices}

\author{Alexey Bufetov}
\address[Alexey Bufetov]{Institute of Mathematics, Leipzig University, Germany}
\email{alexey.bufetov@gmail.com}

\author{Evgeny Obukhov}
\address[Evgeny Obukhov]{}
\email{evgobmm@gmail.com}

\begin{document}
	
	\begin{abstract}		 
		In this note we announce the limiting shape for the height function of a uniformly distributed alternating sign matrix (equivalently, a six-vertex model with equal weights and domain-wall boundary conditions)
	\end{abstract}
	
	\maketitle

\section{Result}		

Let $A^{(N)}=(A_{ij})_{1\leq i,j\leq N}$ be uniform on the alternating sign matrices of size $N$ (see, e.g., \cite{Bressoud99,Propp01}, for definitions). Define
\begin{equation}\label{eq:cornerdef}
	C_N(i,j)=\sum_{a=1}^{i}\sum_{b=1}^{j}A_{ab}, \qquad 0 \le i,j \le N.
\end{equation}
Define the function $c_N: [0,1]^2 \to \R$ via bilinear interpolation of the values $c_N (i/N,j/N):= N^{-1}C_N(i,j)$, $0\leq i,j\leq N$. 

\begin{theorem}\label{thm:main}
	
	There exists a function $c: [0,1]^2 \to \R$ such that
	for every $\varepsilon>0$, there are $a_\varepsilon>0$ and $N_\varepsilon<\infty$ such that
	\begin{equation}\label{eq:LLN}
		\mathrm{Prob} \bigl(\sup_{(x,y) \in [0,1]^2 } \left| c_N (x,y)-c(x,y) \right| >\varepsilon\bigr)
		\leq \exp(-a_\varepsilon N^2),\qquad N\geq N_\varepsilon.
	\end{equation}
	
	The function $c(x,y)$ can be defined in two equivalent ways, given below. 
\end{theorem}

For defining the answer, first recall the arctic curve function
\begin{equation}\label{eq:arctic-quadratic}
	\mathcal Q(x,y):=x^2+y^2-xy-x-y+\frac14.
\end{equation}
We define the function $c(x,y)$ first for a known frozen region:  
\begin{equation}\label{eq:c-quarter-first}
	\begin{aligned}
		c(x,y):=0, \qquad
		0\le x,y<\frac12, \quad \mathcal Q(x,y)\ge0,\\
		c\left(x,\frac12\right)=\frac{x}{2}, \ 
		0\le x\le\frac12,\qquad 
		c\left(\frac12,y\right)=\frac{y}{2}, \ 
		0\le y\le\frac12. 
	\end{aligned}
\end{equation}
Also, $c(x,y)$ on $[0,1]^2$ has the following symmetries:
\begin{equation}\label{eq:c-reflections}
	c(1-x,y)=y-c(x,y),\qquad
	c(x,1-y)=x-c(x,y).
\end{equation}
Thus, it remains to define $c(x,y)$ on $\{ (x,y): 0 < x,y<\frac12,\ \mathcal Q(x,y)<0 \}$.

\section{First description}\label{sec:first-description}

For $B>0$, let $K_B$ be the integral operator on $L^2([-B,B])$ defined via
\begin{equation}\label{eq:fredholm-kernel}
(K_Bf)(t)=\int_{-B}^{B}k(t-s)f(s)\,ds, \qquad \mbox{where $k(z):=\frac{\sqrt3}{2\pi\bigl(1+2\cosh z\bigr)}$}
\end{equation}

We also introduce the function
\begin{equation}\label{eq:ell-def}
    \ell(z):=
    \frac{\sqrt3}
    {4\pi\cosh(z/2)\bigl(2\cosh(z/2)+\sqrt3\bigr)},
    \qquad |\operatorname{Im}z|<\pi.
\end{equation}
For $|v|<\pi$, denote by $Q_{B,v}$ the rank-one operator on
$L^2([-B,B])$ given by
\begin{equation}\label{eq:rank-one}
    (Q_{B,v}f)(t)
    :=\int_{-B}^{B}\ell(s+iv)f(s)\,ds.
\end{equation}
Note that the integral kernel of $Q_{B,v}$ is independent of $t$. Define
\begin{equation}\label{eq:fredholm-dets}
    \tau(B):=\det\bigl(I-K_B\bigr),\qquad
    \tau(B,v):=\det\bigl(I-K_B+Q_{B,v}\bigr),
\end{equation}
and the real-valued function
\begin{equation}\label{eq:P-fredholm}
    \mathcal P(B,v):=\frac{\tau(B,v)}{\tau(B)}-1,
    \qquad B>0,\quad |v|<\pi.
\end{equation}

For $0<x,y<1/2$ with $\mathcal Q(x,y)<0$, define
\begin{equation}\label{eq:phase-first}
    \Phi_{x,y}(B,\theta)
    :=x\mathcal P(B,\theta-\frac{2\pi}{3})
      +y\mathcal P(B,\theta+\frac{2\pi}{3})
      +(x+y-1)\mathcal P(B,\theta),
\end{equation}
for $B>0$ and $|\theta|<\pi/3$. Let
$(B(x,y),\theta(x,y))$ be the unique solution of
\begin{equation}\label{eq:stationary-first}
    \partial_B\Phi_{x,y}(B,\theta)=0,
    \qquad
    \partial_\theta\Phi_{x,y}(B,\theta)=0,
    \qquad B>0,\quad |\theta|<\pi/3.
\end{equation}
We then set
\begin{equation}\label{eq:c-liquid-first}
    c(x,y):=
    \Phi_{x,y}\bigl(B(x,y),\theta(x,y)\bigr).
\end{equation}

This defines $c$ in the liquid part of the lower-left quarter of the
square.

\section{Second description}\label{sec:second-description}

We now give an equivalent description in terms of a solution
of an equivalent form of Painlev\'e VI.  The resulting
answer is closely related to the general
framework of Kenyon-Prause for gradient variational problems and limit shapes
\cite{KenyonPrause22,KenyonPrause24}.

Let $r(B)$ be the unique real-analytic solution for $B\geq0$ of
\begin{equation}\label{eq:r-pvi}
    r''
    =r-2r^3
      +2r\coth(3B)\sqrt{(r')^2+r^4-r^2},
    \qquad B>0,
\end{equation}
where the square-root is nonnegative, with
\begin{equation}\label{eq:r-pvi-data}
    r(0)=0,
    \qquad
    r'(0)=\frac{\sqrt3}{\pi}.
\end{equation}
Define
\begin{equation}\label{eq:rho-pvi}
    \rho(B):=\frac{r(B)}{\sinh(3B)},
    \qquad
    \rho(0)=\frac{1}{\pi\sqrt3}.
\end{equation}

Let
\[
    \mathbf F(B,\theta)
    =\bigl(F_0(B,\theta),F_1(B,\theta),F_2(B,\theta)\bigr)
\]
be the unique real-analytic solution for $B\geq0$ and
$|\theta|<\pi/3$ of
\begin{equation}\label{eq:F-pde}
    \partial_B^2\mathbf F+\partial_\theta^2\mathbf F
    =2\rho(B)\,\partial_B\mathbf F,
\end{equation}
with Cauchy data
\begin{equation}\label{eq:F-data}
    \mathbf F(0,\theta)=0,
    \qquad
    \partial_B\mathbf F(0,\theta)
    =\frac{
        \bigl(1,\cos(\theta/2),\sin(\theta/2)\bigr)
      }{\cos(3\theta/2)}.
\end{equation}

For $(x,y)\in[0,1]^2$, set
\begin{equation}\label{eq:XY-coordinates}
    X:=\frac{2-x-y}{\sqrt3},
    \qquad
    Y:=x-y.
\end{equation}
Then
\begin{equation}\label{eq:delta-circle}
    \mathcal Q(x,y)=\frac34\bigl(X^2+Y^2-1\bigr).
\end{equation}
For $0<x,y<1/2$ with $\mathcal Q(x,y)<0$, define
\begin{equation}\label{eq:phase-second}
    \Psi_{x,y}(B,\theta)
    :=F_0(B,\theta)-X F_1(B,\theta)-Y F_2(B,\theta).
\end{equation}
Let $(B(x,y),\theta(x,y))$ be the unique solution in
$B>0$, $|\theta|<\pi/3$ of
\begin{equation}\label{eq:stationary-second}
    \partial_B\Psi_{x,y}(B,\theta)=0,
    \qquad
    \partial_\theta\Psi_{x,y}(B,\theta)=0.
\end{equation}
Then
\begin{equation}\label{eq:c-liquid-second}
    c(x,y)
    :=\frac{3}{2\pi}\,
      \Psi_{x,y}\bigl(B(x,y),\theta(x,y)\bigr).
\end{equation}

The equivalence of the two descriptions follows from an explicit linear
relation between $\bigl(F_0,F_1,F_2\bigr)$ and the three shifted functions
$\mathcal P(B,\theta)$ and
$\mathcal P\left(B,\theta\pm\frac{2\pi}{3}\right)$.

\section{Brief History}

Alternating sign matrices were introduced by Mills, Robbins, and Rumsey in the early 1980s \cite{MillsRobbinsRumsey83}.  They conjectured the product formula
\[
    \prod_{j=0}^{N-1}\frac{(3j+1)!}{(N+j)!}
\]
for the number of $N\times N$ ASMs.  The conjecture was proved by Zeilberger \cite{Zeilberger96}; Kuperberg soon gave a second proof based on the correspondence with square ice (the six-vertex model with domain-wall boundary conditions) \cite{Kuperberg96}.  The enumeration story and its surrounding conjectures are described in Bressoud's book \cite{Bressoud99}; see also Propp's survey \cite{Propp01} for the many equivalent realizations of ASMs.

Numerical sampling of the uniformly random ASM gave strong evidence for macroscopic phase separation, see Allison-Reshetikhin \cite{AllisonReshetikhin05} and Wieland \cite{WielandASM}.  Colomo and Pronko derived a conjectural arctic curve for uniformly random ASMs from the emptiness-formation probability and a condensation hypothesis \cite{ColomoPronkoLimit10}, and subsequently obtained the general disordered-regime six-vertex formula \cite{ColomoPronkoArctic10}.  Colomo and Sportiello introduced the tangent method, which rederived the square-domain curve and predicted arctic curves on substantially more general domains \cite{ColomoSportiello16}.  Aggarwal proved the arctic curve phenomenon for uniformly random ASMs via rigorously justifying the tangent-method mechanism in this setting \cite{Aggarwal20}. 

A broader surface-tension variational framework for six-vertex limit shapes with fixed boundary conditions was formulated and studied by Palamarchuk and Reshetikhin \cite{PalamarchukReshetikhin08}.  The integrable structure of the corresponding limit-shape equations was further developed by Reshetikhin and Sridhar \cite{ReshetikhinSridhar17}, and, for inhomogeneous six-vertex models, by Keating, Reshetikhin, and Sridhar \cite{KeatingReshetikhinSridhar22}. This general variational framework seems to be intricate, and it is not easy to make explicit computations via it for our case.

\section{Numerical checks}

The following table compares the values of $c(x,y)$ obtained independently from the first and second descriptions with numerical data from Wieland \cite{WielandASM}. The last two columns are sample averages for $N=599$ (300 samples) and $N=999$ (100 samples), respectively.

\begin{center}
\scriptsize
\setlength{\tabcolsep}{3.5pt}
\renewcommand{\arraystretch}{1.12}
\begin{tabular}{c|c|c|c|c}
$(x,y)$ & First description & Second description & $N=599$ & $N=999$ \\
\hline
$(0.10,0.20)$ & $0.001527589483$ & $0.001527589483$ & $0.001522381226$ & $0.001572572573$ \\
$(0.10,0.30)$ & $0.014217459146$ & $0.014217459146$ & $0.014170830984$ & $0.014167067067$ \\
$(0.15,0.35)$ & $0.039892658925$ & $0.039892658925$ & $0.039978029157$ & $0.039987537538$ \\
$(0.20,0.30)$ & $0.044868991022$ & $0.044868991022$ & $0.044914657093$ & $0.044830230230$ \\
$(0.25,0.25)$ & $0.046503273995$ & $0.046503273995$ & $0.046497506993$ & $0.046393268268$ \\
$(0.20,0.40)$ & $0.071972768004$ & $0.071972768004$ & $0.071921653145$ & $0.071910310310$ \\
$(0.25,0.40)$ & $0.092937048974$ & $0.092937048974$ & $0.092862407245$ & $0.092892392392$ \\
$(0.30,0.40)$ & $0.114131348300$ & $0.114131348300$ & $0.114043472677$ & $0.114113513514$ \\
$(0.35,0.40)$ & $0.135481175841$ & $0.135481175841$ & $0.135392934713$ & $0.135523423423$ \\
$(0.40,0.40)$ & $0.156934672482$ & $0.156934672482$ & $0.157011965316$ & $0.156934134134$ \\
\end{tabular}
\end{center}

\subsection*{Comments}

The proof of the result was assisted by ChatGPT. 

It starts from the general variational approach to random surfaces developed by Sheffield \cite{Sheffield05} and, in the form used here, by Lammers-Tassy \cite{LammersTassy24}. The other important ingredients include ideas and results from Alessandrini \cite{Alessandrini87}, Duminil-Copin-Kozlowski-Krachun-Manolescu-Tikhonovskaia \cite{DCKKMT22}, Kitanine-Maillet-Terras \cite{KMT99}, and Kenyon-Prause \cite{KenyonPrause22}, \cite{KenyonPrause24}. It also uses the already known arctic curve, conjectured by Colomo-Pronko \cite{ColomoPronkoLimit10,ColomoPronkoArctic10} and proved by Aggarwal \cite{Aggarwal20}, so it does not produce an independent proof of that curve.

The proof has been verified in Lean modulo several published results. The Lean formalization will be made public together with the detailed version of the paper. 

\subsection*{Acknowledgements}

A.~Bufetov was partially supported by the European Research Council (ERC), Grant Agreement No. 101041499.

\end{document}